\documentclass{llncs}
\usepackage{llncsdoc}

\usepackage{amsmath,amssymb,bbm}
\usepackage{mathrsfs}
\usepackage{graphicx}

\def\o{\circ\,} 
 
\def\al{\alpha} 
\def\be{\beta} 
\def\ga{\gamma} 
\def\de{\delta} 
\def\ep{\varepsilon}

\def\si{\sigma} 
 
\def\ph{\varphi}

\def\De{\Delta}

\def\i{^{-1}}

\def\p{\partial}

\let\on=\operatorname

\def\R{\mathbb{R}}
\let\wh=\widehat

\begin{document}

\title{A New Riemannian Setting for Surface Registration}
\author{Martin Bauer\inst{1}, Martins Bruveris\inst{2}}
\institute{Fakult\"at f\"ur Mathematik, Universit\"at Wien, 
Nordbergstrasse 15, A-1090 Wien
\and Dep. of Mathematics, Imperial College, London SW7~2AZ, UK}

\maketitle

\begin{abstract}
We present a new approach for matching regular surfaces in a Riemannian setting. 
We use a Sobolev type metric on deformation vector fields which form the tangent bundle to the space of surfaces. In this article we compare our approach with the diffeomorphic matching framework. In the latter
approach a deformation is prescribed on the ambient space, which then drags along an embedded surface. In contrast our metric is defined directly on the deformation vector field and can therefore be called an {\it inner metric}. We also show how to discretize the corresponding geodesic equation and compute the gradient of the cost functional using finite elements.

\end{abstract}

\section{Introduction}

The field of computational anatomy concerns itself with the study and classification of the variability of biological shapes, including their statistical variance. Since the space of all shapes is inherently nonlinear, the usual methods of linear statistics cannot be applied. In particular, the addition of two surfaces cannot be meaningfully defined. One way to overcome this difficulty is to introduce a Riemannian structure on the space of shapes, which locally linearizes the space and allows the development of statistical methods that are analogous to the linear case. This approach was taken, e.g., in \cite{Fletcher2008}. In the Riemannian setting, the average of two shapes may be defined as the middle point of a geodesic joining these two shapes. In a similar way one may define the corresponding geodesic mean of a collection of $n$ shapes.

One class of shapes which are of interest in computational anatomy consists the surfaces embedded in $\R^3$. The cortical surface, the surfaces of hippocampi, thalami, and nasal cavities are all examples of shapes which are represented as two dimensional surfaces in $\R^3$. This is also an example, where the Riemannian setting may be applied to study collections of shapes.

Throughout the last decade various (Riemannian) metrics have been proposed. They include  a metric that preservers local rigidity \cite{Kilian2007}, a generalization of the elastic metric for curves to higher dimensions \cite{Mio2005,Liu2010}, a metric inspired by a continuum mechanics which is defined in the interior of a two dimensional shape \cite{Wirth2011,Fuchs2009} and a metric based on the square-root representation of surfaces \cite{Kurtek2010}.
Other approaches include the representation of surfaces via densities \cite{Peter2008} and metrics defined on surfaces via level sets \cite{Bertalmio2001,Memoli2004}.

Another method for comparing anatomical shapes in the Riemannian setting is the method of large deformation diffeomorphic metric matching (LDDMM), based on the deformable template paradigm of Grenander 
\cite{Grenander1993}. In this setting, a template shape is matched to a target shape by finding a transformation in a suitable group of deformations of the ambient space that transforms the template into the target. This approach has been systematically developed in \cite{Beg2005,Bruveris2011,Glaunes2008,Risser2010,Vialard2011} and applied to various problems in computational anatomy. Registering two surfaces in this framework involves finding a diffeomorphism of the whole space, which transforms one surface to the other. Because of its widespread use in the field of computational anatomy for registering volumetric images we will use LDDMM as a reference to highlight the features of our proposed framework.  

In this paper we propose a different way of defining a Riemannian structure on the space of surfaces, which also provides the full range of tools for nonlinear statistics. We use a Sobolev type norm to enforce regularity of the deformation vector field and measure the cost of the deformation. Our approach to characterizing a deformation is intrinsic to the surface, rather than resulting from a transformation of the surface induced by a deformation of the ambient space in which the image is embedded. For this reason, the Riemannian metrics used here are called {\it inner metrics} as opposed to the {\it outer metrics} used in LDDMM, where deformations are imposed via the ambient space. Other examples for inner metrics can be found in \cite{Mio2005,Kurtek2010,Wirth2011,Fuchs2009,Michor118}.

Inner metrics of Sobolev type on planar curves were introduced and studied previously in \cite{Michor107,Michor111}. Recently they were generalized to surfaces and higher dimensional hypersurfaces in Euclidean space in \cite{Michor119}. The numerical implementation of matching with these metrics  differs from LDDMM, because the metric on the tangent space at each surface depends nonlinearly on the surface. This means the metric will change adaptively as one moves around in shape space. This adaptive property is in marked contrast to LDDMM, where the metric is defined on the diffeomorphisms and projected down to the shapes, so it doesn't depend on the particular shape.

The outline of the paper is as follows. In Sec. \ref{math} we review the registration problem for surfaces and recall how it is solved using outer metrics in LDDMM. Then we present the approach via inner metrics of Sobolev type and point out the differences between the two methods. For definiteness, we will concentrate our attention on the Sobolev metric of order one. In Sec. \ref{disc} we discuss how to discretize and implement the geodesic equations for this metric and how to solve the registration problem via geodesic shooting. Finally, in Sec. \ref{example} we show how this metric performs in some examples using synthetic data.

\section{The Mathematical Formulation}
\label{math}

We are dealing with the registration of parametrized regular surfaces. Such a surface is given by a smooth function $q: M \to \R^3$ from a model surface $M$ into the Euclidean space $\R^3$. We will consider different choices of the model surface $M$ in this paper: the plane sheet $M=[0,1]\times[0,1]$, the cylinder $M=S^1 \times [0,1]$ and the torus $M=S^1\times S^1$. Another interesting choice would be the sphere $M=S^2$, which however is not considered in this paper. The metric can be defined in the same way as for the other topologies, however the numerical treatment is more challenging, because the sphere cannot be covered by a single global coordinate chart. We require the parametrization of the surface $q$ to be regular in the following sense: at each point $x \in M$ the partial derivatives $\frac {\p q}{\p x^1}$, $\frac {\p q}{\p x^2}$ are required to be linearly independent. We will denote the space of all such surfaces by $\mathscr S$.

\subsection{Registration with LDDMM}
\label{sec_lddmm}

In the LDDMM framework, the registration of a template surface $q_0$ to a target surface $q_{\mbox{targ}}$ involves finding a curve $\ph_t$ of diffeomorphisms of the ambient space $\R^3$, whose deformation carries the template surface to the target surface. Mathematically, one constructs these deformations using time dependent vector fields $v_t(y)$, which generate $\ph_t$ as their flow, i.e.
\begin{equation}
\p_t \ph_t = v_t(\ph_t) \enspace.
\end{equation}
The registration problem consists of finding a vector field, which minimizes the following sum of a kinetic energy and a matching term
\begin{equation}
\label{en_lddmm}
E(v_t) = \frac 12 \int_0^1 \| v_t \|_V^2\,dt + \frac{1}{2\si^2} d(\ph_1(q_0), q_{\mbox{targ}}) \enspace.
\end{equation}
The kinetic energy is usually measured using a norm $\|.\|_V$ defined on a reproducing kernel Hilbert space $V$ of vector fields on $\R^3$ with kernel $K$. The norm is then given by $\| u \|^2_V = \int u \cdot K\i \star u\,dx$. We will discuss possible choices of the matching term in Section \ref{matching}.

It is possible to reduce the complexity of the problem, since one can show that minimizing vector fields $v_t$ have to obey an evolution equation. This enables us 
to describe the whole vector field $v_t$ by knowing only its value $v_0$ at time $t=0$. The equations
\begin{subequations}
\label{geod_lddmm}
\begin{align}
\p_t q_t(x) &= v_t(q_t(x)) \label{drag} \\
\p_t p_t(x) &= (Dv_t(q_t(x)))^T p_t(x) \label{drag_mom} \\
v_t(y) &= \int_M p_t(x)K(y-q_t(x))\,dx \label{conv_mom}
\end{align}
\end{subequations}
are given in terms of a momentum $p_t$, which lives on the surface. The momentum is convolved with the kernel in \eqref{conv_mom} to reconstruct the minimizing vector field, which defines the deformation of the ambient and drags the surface along in \eqref{drag}. Details about LDDMM and the reformulation as evolution equations can be found in \cite{Beg2005,Younes2009}.

\subsection{Registration with Inner Metrics}

\begin{figure}
\centering
\includegraphics[height=4cm]{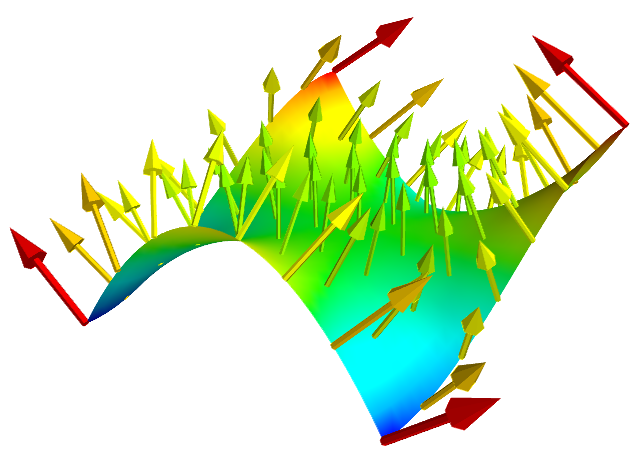}
\hfill
\includegraphics[height=4cm]{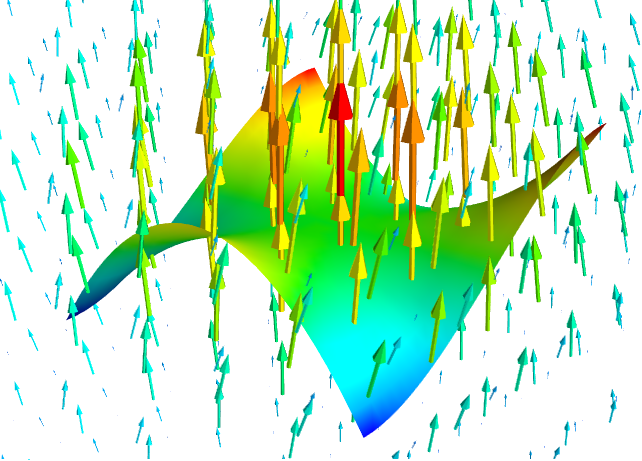}
\caption{For inner metrics the vector field governing the deformation is defined directly on the surface (left picture). In contrast the LDDMM approach defines a deformation vector field on all of $\R^3$ (right picture). The latter vector field deforms the ambient space and in the process induces a deformation of the surface.}
\label{vector_fields}
\end{figure}

We propose to use a different approach, described from a mathematical point of view in \cite{Michor118,Michor119}. In this approach we describe the deformation of the surface directly, without assuming an underlying deformation of the whole space. In our approach we will replace \eqref{drag} by
\begin{equation}
\p_t q_t(x)= u_t(x) \enspace,
\end{equation}
where $u_t(x) \in C^\infty(M, \R^3)$ is a time dependent vector field, defined only on the surface. Note the difference between the vector field $v_t(y)$ which is defined on $\R^3$ and $u_t(x)$, which is defined on the model space $M$, 
c.f. Fig.~\ref{vector_fields}.

In the new framework the registration problem still consists of minimizing the following sum of a kinetic energy term and a matching functional
\begin{equation}
\label{en_inner}
E(u_t) = \frac 12 \int_0^1 \langle u_t, u_t \rangle_{q_t}\,dt + \frac{1}{2\si^2} d(q_1,q_{\mbox{targ}}) \enspace .
\end{equation}
The kinetic energy is measured via an inner product $\langle .,. \rangle_{q_t}$ on the space of vector fields along the surface. This inner product can and usually will depend nonlinearly on the surface $q_t$ itself. This is another difference with the LDDMM framework, where the kinetic energy of the vector fields didn't depend on the surfaces that were matched. There is a whole variety of inner metrics, one can choose from. We will concentrate in this paper on an $H^1$-type metric, which will be discussed in Sec.\ref{h1metric}.
Again we postpone the discussion of the matching term to Sec. \ref{matching}.

\subsection{Geometry}

Both the LDDMM approach and the inner metrics can be seen as a special case of a constructions in the 
general framework of Riemannian geometry. In the case of inner metrics the collection of inner products $\langle .,. \rangle_{q}$ defines a Riemannian metric on the space $\mathscr S$ of all surfaces. In Riemannian geometry, curves, which minimize the energy $\int_0^1 \langle u_t, u_t \rangle_{q_t}\,dt$ for fixed endpoints $q_0, q_1$ are called geodesics. We can see that minimizers of \eqref{en_inner} have to be geodesics in the space $\mathscr S$.

In LDDMM the inner product on the space of vector fields $V$ also defines a Riemannian metric, this time on the group of diffeomorphisms. Therefore the minima of the registration problem \eqref{en_lddmm} generates geodesics in the diffeomorphism group.

Why is it advantageous to work in a Riemannian setting? In this setting, the minima of the matching energy are geodesics, so one may describe the nonlinear space of shapes $\mathscr S$ in terms of the initial velocity or momentum of a geodesic, which is an element in a linear vector space. This is possible, because geodesics obey an evolution equation like \eqref{geod_lddmm}. Using the initial velocity, which encodes the whole solution of the registration problem, we are able to view the space of surfaces from the template surface $q_0$ as a linear space. This enables us to use statistics, compute average surfaces and measure distances.

\subsection{$H^1$-type Metric on Surfaces}\label{h1metric}

We will consider a metric on the space of surfaces, which is the analogue of the $H^1$-norm $\|f\|_{H^1}^2 = \int_{\R^2} |f(x)|^2 + |\nabla f(x)|^2\, dx$ for functions on $\R^2$. We will replace functions on $\R^2$ by vector fields living on the curved surface $q$ and adjust the definition of the $H^1$-norm to take into account the curved nature of the surface $q$. Since $q$ is a surface in $\R^3$, we can measure angles and distances of vectors tangent to the surface, using the Euclidean inner product on $\R^3$. At each point $q(x)$ of the surface we also have a canonical basis for the plane tangent to $q$, given by the vectors $\frac {\p q}{\p x^1}$, $\frac {\p q}{\p x^2}$. We  denote 
the inner product induced on the tangent plane to the surface by $g$.
This inner product has the following coordinate matrix with respect to the basis 
$\frac {\p q}{\p x^1}$, $\frac {\p q}{\p x^2}$:
\begin{align}
(g_{ij}) &= \left( \sum_{k=1}^3 \frac {\p q^k}{\p x^i} \frac {\p q^k}{\p x^j} \right)_{i,j=1}^2, 
& (g^{ij}) &= (g_{ij})\i, & \on{vol}(g) &= \sqrt{\det(g_{ij})}\enspace.
\end{align}
We denote by $(g^{ij})$ the inverse matrix of $(g_{ij})$ and by $\on{vol}(g)$ the volume density of the surface $q$ with respect to the metric $g$. 

In this paper we will use the $H^1$-type inner metric on surfaces defined by the expression
\begin{align}
\label{h1_metric}
\langle u,v\rangle_{q}:&=\sum_{k=1}^3 \int_M u^kv^k
+\al^2 g\left(\on{grad}^g(u^k),\on{grad}^g(v^k)\right)\on{vol}(g)\,dx \\
&= \sum_{k=1}^3 \int_M u^kv^k + \al^2 g^{ij} \frac {\p u^k}{\p x^i} \frac {\p v^k}{\p x^j} \on{vol}(g) \,dx \enspace.
\end{align}
 One reason to use this generalization of the $H^1$-metric is, that this metric is invariant under reparametrizations of the surface and only depends on the image $q(M)$ as a subset of $\R^3$, in a similar way as the length of a curve in two dimensions only depends on the image of the curve, and not on a particular parametrization. This is necessary, if one wants to match unparametrized surfaces. For this task it is possible to use the same framework with this metric, only the matching term has to be chosen to beinvariant under reparametrizations.

The constant $\al$, which appears in the metric, is a parameter, which has to be chosen for each problem. It represents the characteristic length scale, at which deformations take place. Another interpretation of $\al$ is the scale across which the momentum is smoothed, when passing from momenta to velocities. It can be compared with the kernel size in LDDMM.

Other choices for the metric are possible. One could use Sobolev type metrics of higher order as in \cite{Michor119} or multiply the components of the metric with a function depending on geometric quantities of the surface, like area, mean or Gaussian curvature as was done in \cite{Michor118}.

\subsection{The Matching Functional}
\label{matching}

There are different possible choices for the matching term. In this paper we will use the squared $L^2$-distance
\begin{equation}
d(q_0, q_1) = \int_M |q_0(x) - q_1(x) |^2\,dx.
\end{equation}
Since we are dealing with parametrized surfaces this is a natural choice for the matching functional.

When matching unparametrized surfaces, natural choices of the matching functional would include currents, see \cite{Glaunes2008}, or one could use the reparametrization framework of \cite{Cotter2010}.

\section{Discretization}
\label{disc}

In this section we will describe how to discretize the variational problem \eqref{en_inner} and compute the optimal path between two surfaces. Starting with an initial guess for the velocity $u_0^1$, we will use a gradient descent scheme 
\[ u_0^{i+1} = u_0^i - \ep \nabla_{u_0} E (u_0^i) \]
to converge towards the initial velocity of the optimal geodesic. The discretization thus consists of two parts:
\begin{itemize}
\item compute the geodesic, given the initial velocity to evaluate $E(u_0^i)$
\item compute the gradient $\nabla_{u_0} E (u_0^i)$ to update the initial velocity.
\end{itemize}

We show how to discretize the geodesic equation in Sec. \ref{geod_sec} and how to compute the gradient in Sec. \ref{grad_sec}.

\subsection{The Geodesic Equation}
\label{geod_sec}

We discretize the time-evolution of the surface $q(t)$ using the explicit Euler method
\begin{equation}
  \label{discr_geod}
  q_{i+1} = q_i + \De t u_i \enspace,
\end{equation}
where $q_i = q(i\De t)$ is the discretized version of the curve and $\De t = 1/N$ is the time step, if we divide the interval $[0,1]$ into $N$ parts. To compute $u_i$ we note that a geodesic is a critical point of the energy $E(u_i) = \frac{\De t}{2} \sum_{i=0}^{N-1} \langle u_i, u_i \rangle_{q_i}$, i.e. $\nabla_{u_i} E(u_i) =0$. Following \cite{LeDimet1997} we introduce the Lagrangian multiplier $p_i$ in the discrete variational principle
\begin{equation}
  E(u_0,\ldots,u_{N-1}) = \sum_{i=0}^{N-1} \frac{\De t}{2} \langle u_i, u_i \rangle_{q_i} + \langle p_i, q_{i+1} - q_i - \De t u_i \rangle_{L^2} \enspace.
\end{equation}
and take variations. From variations in $u_i$  we see that $p_i$ is the momentum dual to the velocity $u_i$ in the sense that 
$\langle u_i, \de u_i \rangle_{q_i} = \langle p_i, \delta u_i \rangle_{L^2}$
and we obtain the evolution equation for $u_i$ in the form
\begin{equation}
  \label{u_impl}
  \langle u_{i+1}, \de q_{i+1} \rangle_{q_{i+1}} = \langle u_i, \de q_{i+1} \rangle_{q_i} + \De t \left\langle \frac{\de \ell}{\de q}(u_{i+1}, u_{i+1}; q_{i+1}), \de q_{i+1} \right\rangle_{L^2} \enspace,
\end{equation}
with $\de q_{i+1}$ an arbitrary variation. Here we use the notation
\begin{equation}
  \ell(u,v;q) = \frac12 \langle  u, v \rangle_q \enspace.
\end{equation}
We denote by $\frac{\de \ell}{\de q}$ the variational derivative of $\ell(u,v;q)$ with respect to the variable $q$, defined via 
\begin{equation}
  \left\langle \frac{\de \ell}{\de q}(u, v; q), \de q \right\rangle_{L^2} = \lim_{h\to 0} \frac{\ell(u,v;q+h\de q)-\ell(u,v;q))}{h} \enspace.
\end{equation}

Equation \eqref{u_impl} is an implicit time step for $u_{i+1}$, since $u_{i+1}$ appears on the right hand side in a quadratic term. To make computations easier and avoid having to solve a nonlinear equation, we changed the right hand side to an explicit Euler time step
\begin{equation}
    \label{discr_geod2}
 \langle u_{i+1}, \de q \rangle_{q_{i+1}} = \langle u_i, \de q \rangle_{q_i} + \De t \left\langle \frac{\de \ell}{\de q}(u_{i}, u_{i}; q_{i}), \de q \right\rangle_{L^2} \enspace.
\end{equation}

\subsection{Computing the Gradient}
\label{grad_sec}

Given \eqref{discr_geod}, \eqref{discr_geod2} for the evolution of a geodesic we again use method of adjoint equations from \cite{LeDimet1997} to compute the gradient of the energy with respect to the initial velocity. The resulting equations for the variables $\wh u_i, \wh v_i$ have to be integrated backwards in time
\begin{equation}
\label{adj_eq}
\begin{split}
\langle \wh u_{i}, \de u_i \rangle_{q_i} &=
\langle \wh u_{i+1}, \de u_i \rangle_{q_i} + \De t \langle \wh v_{i+1}, \de u_i \rangle_{q_{i+1}} + 2 \De t \langle \frac{\de\ell}{\de q}(u_i,\de u_i; q_i),\wh u_{i+1} \rangle \\
\langle \wh v_{i}, \de q_i \rangle_{q_{i}} &=
\langle \wh v_{i+1}, \de q_i \rangle_{q_{i+1}} + 2 \left\langle \frac{\de\ell}{\de q}(\wh u_{i+1} - \wh u_{i},u_i; q_i),\de q_i \right\rangle \\
&\phantom{\mathrel{=}}+\De t \left\langle \frac{\de^2\ell}{\de q^2}(u_i,u_i; q_i),(\wh u_{i+1}, \de q_i) \right\rangle
\end{split}
\end{equation}
with the initial conditions
\begin{align}
  \wh u_{N} &= 0 &
  \langle \wh v_{N}, \de q_N \rangle_{q_{N}} &= -\frac 1{\si^2} \langle q_N - q_{\mbox{targ}}, \de q_N \rangle
\end{align}
at time $t=1$. The gradient is then given by
\begin{equation}
  \nabla_{u_0} E(u_0) = u_0 - \wh u_0 \enspace.
\end{equation}




\section{Numerical Experiments}
\label{example}

We implemented the geodesic and adjoint equations \eqref{discr_geod}, \eqref{discr_geod2} and \eqref{adj_eq} in Python using the finite element library FEniCS \cite{Fenics2010}. All model manifolds ($[0,1]\times[0,1]$, $S^1 \times [0,1]$, $S^1\times S^1$) were modelled on the rectangle $[0,1]\times[0,1]$ with periodic boundary conditions prescribed where neccessary. The domain was subdivided into a regular triangular mesh, on which Lagrangian finite elements of order 1 were defined.

\begin{figure}
\centering
\includegraphics[height=4.5cm]{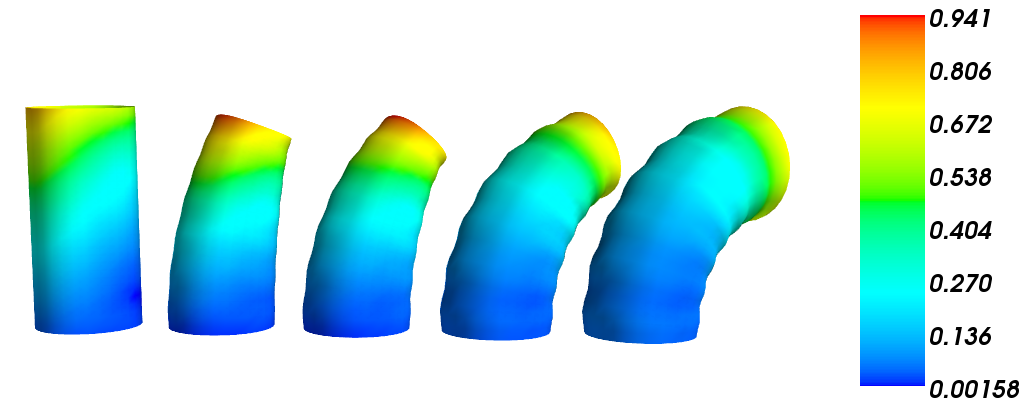}
\caption{Samples are shown from a geodesic in the space of surfaces between a straight cylinder and a bent cylinder with ripples, at time points $t=0, 0.3, 0.5, 0.8, 1$. The color encodes the Euclidean length of the deformation vector field at each point of the surface.}
\label{ex1}
\end{figure}

In the first example we apply our method to compute the geodesic path between two shapes, which includes both large and small deformations. The template shape is a straight cylinder of height 1 and radius 0.25, which is discretized using a regular triangular mesh of $2\times30\times30$ elements. The target shape is a cylinder, which is bent by $90^\circ$ and has 5 small ripples added to it along the vertical axis. Compared to the bending the ripples constitute a small and local deformation of the shape. The target shape is discretized in the same way as the template. We use $\alpha=0.6$ as the length scale parameter and 10 time steps for the time integration. The gradient descent takes 80 steps to converge to an $L^2$-error of $0.008$. We can see in Fig. \ref{ex1} that both the large and the small deformations are captured by the geodesic.

\begin{figure}
\centering
\includegraphics[height=6cm]{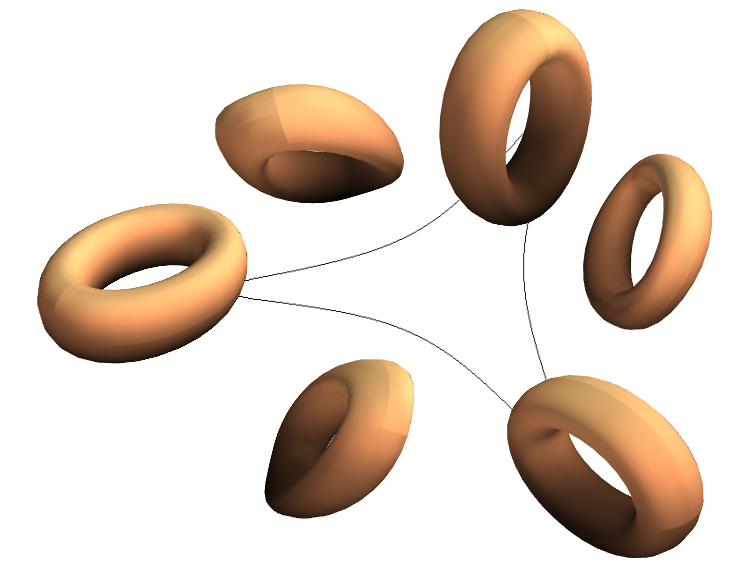}
\caption{This figure shows a geodesic triangle in the space of surfaces with asymmetric tori as vertices. The tori along the edges are the middle points of the geodesics connecting the vertices. One can see that the shapes tend to shrink along the geodesics before expanding again towards the ends. This effect implies negative curvature in this region of shape space.}
\label{ex2}
\end{figure}

In the second example we want to illustrate the curved nature of shape space. To do so we pick three asymmetric tori, lying in different positions in space. Each two tori differ by a composition of two rigid rotations. We compute the geodesics between each pair of tori to measure the angles and side lengths of the triangle with the tori as vertices and the geodesics as edges. By comparing the sum of the angles with $\pi$ one can estimate, whether the curvature of shape space along the plane containing the triangle is positive or negative. We measured $\al = 33.766^\o$, $\be = 34.802^\o$ and $\ga = 34.675^\o$. The sum $\al+\be+\ga = 103.243^\o$ is smaller than $180^\o$, which indicates that the space is negatively curved in this area (c.f.~\cite[Sec.~5.4]{Michor98}). In negatively curved spaces geodesics tend to be attracted towards a common point. In this example the geodesics are attracted towards the surface, which is degenerated to a point. We can see in Fig. \ref{ex2} that the midpoints of the geodesics between the vertices are slightly shrunk. This is another indication for the negatively curved nature of the space.

\begin{figure}
\centering
\includegraphics[height=6cm]{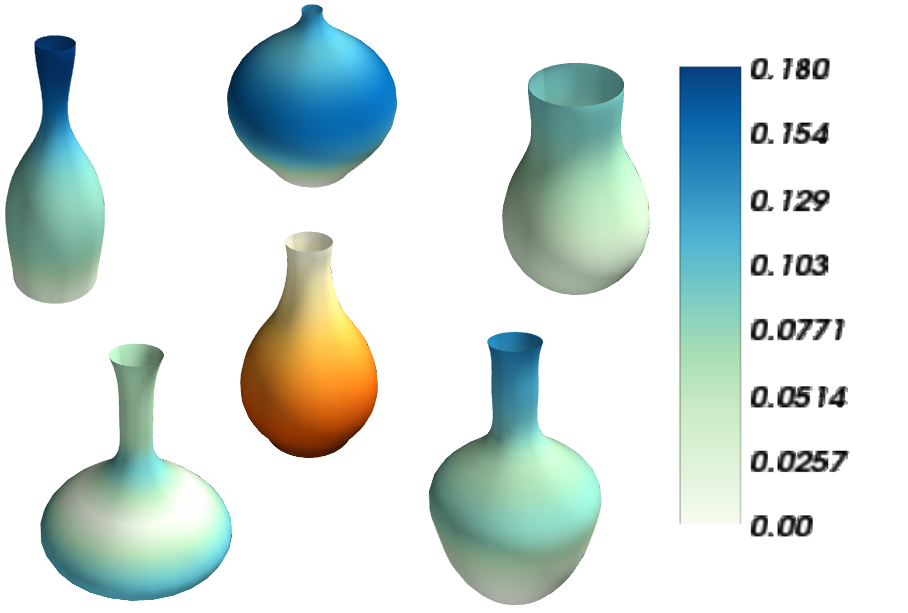}
\caption{In this figure we show the Karcher mean of five vase-shaped objects. The mean shape, which is displayed in the center of the figure is computed using an iterated shooting method. The colored regions  on the averaged shapes encode the Euclidean length of the initial velocity of the geodesic, which connects each shape to the mean. The color of the mean was chosen for artistic purposes only.}
\label{ex3}
\end{figure}

In the third example we show, that our framework is indeed capable to do nonlinear statistics on shape space. We generate five sample shapes and compute the mean shape between them. The five shapes are cylindrical vases with an open top and bottom, discretized again using a triangular mesh of $2\times30\times30$ elements. As the initial guess for the mean we use a straight cylinder. First we register this initial shape to the five target shapes, compute the average of the initial velocities and then shoot with this average velocity to obtain a next guess for the mean shape. We iterate this procedure until the average velocity is close to zero. This method of computing the Karcher mean was proposed in \cite{Fletcher2008}. After four iterations we obtained an average velocity with norm $0.006$. As can be seen in Fig. \ref{ex3} the average shape indeed combines the characteristics of the five shapes.

\section{Conclusions}

In this paper we propose a metric to match regular surfaces in a Riemannian setting. Although this metric has been studied from a mathematical point of view in several papers, including \cite{Michor119,Michor107}, so far it hasn't been applied to problems in computational anatomy. The aim of this work is to argue, that this is a promising approach, which is worthwhile to be studied further. 

In contrast to LDDMM, where surfaces are deformed via a deformation of the ambient space, in this approach the deformation is prescribed directly on the surface, while the ambient space stays constant. Because of this we call this approach matching with {\it inner metrics} as opposed to LDDMM, which can be described as matching with {\it outer metrics}. Other inner metrics, which have been proposed in the literature include \cite{Mio2005,Kurtek2010,Wirth2011,Fuchs2009,Michor118}
We show how to discretize the geodesic equation and how to compute the gradient of the matching functional with respect to the initial velocity. In the last part of the paper we present numerical results on synthetic data of different topologies, which demonstrate the versatility and applicability of our method.

 At the present we applied this metric to match parametrized surfaces, which is an unwelcome restriction in applications. This is not a restriction of the framework itself, but only of the matching functional. By choosing a matching functional, which is independent of the parametrization of the surface, one can apply the same framework to unparametrized surfaces. In future work we plan to implement this capability and test the method on real anatomical data.

\subsection*{Acknowledgments}
We are enormously grateful to our friends and colleagues Colin Cotter, Phillip Harms, Darryl Holm, Peter Michor and David Mumford for their thoughtful comments and encouragement in the course of this work. 

We also acknowledge partial support by the Royal Society of London Wolfson Award, the European Research Council Advanced Grant, the Imperial College London SIF Programme and the Austrian Science Fund (FWF).

\end{document}